\documentclass[12pt]{article}
\usepackage{amsmath}

\usepackage{amsthm} 
\usepackage{amscd}
\theoremstyle{plain} 
\newtheorem{theorem}{Theorem}[section]
\newtheorem{corollary}[theorem]{Corollary}

\newtheorem{lemma}[theorem]{Lemma}
\newtheorem{proposition}[theorem]{Proposition}

\numberwithin{equation}{section}

\newcommand{\al}{\alpha}

\newcommand{\sg}{\sigma}
\newcommand{\br}{\mathbf {R}}

\newcommand{\bz}{\mathbf {Z}}

\newcommand{\ot}{\otimes}

\newcommand{\g} {\mathsf g }

\begin{document}
\title{A Comparison of Leibniz and Cyclic Homologies}
\author{Jerry M. Lodder}
\date{}
\maketitle

\noindent
{\em Mathematical Sciences, Dept. 3MB, Box 30001, New Mexico State
University, Las Cruces NM, 88003, U.S.A. }

\noindent
e-mail:  {\em jlodder@emmy.nmsu.edu}

\bigskip
\noindent
{\bf Abstract.}  We relate Leibniz homology to cyclic homology
by studying a map from a long exact sequence in the Leibniz
theory to the ISB exact sequence in the cyclic theory.  This
provides a setting by which the two theories can be compared
via the 5-lemma.  We then show that the Godbillon-Vey invariant,
as detected by the Leibniz homology of formal vector fields, maps
to the Godbillon-Vey invariant as detected by the cyclic homology
of the universal enveloping algebra of these vector fields.  



\section{Introduction}

Leibniz homology, $HL_*$, was introduced in part as a tool
for studying the obstruction to the periodicity of algebraic
K-theory.  See the introduction to \cite{LodayD}.  As Hochschild
homology, $HH_*$, measures the failure of cyclic homology, $HC_*$, to
be periodic, Leibniz homology should measure the failure of a 
conjectured motivic complex to be periodic.  In this paper we 
begin with a map from Leibniz homology to Hochschild homology,
 and then show that a natural long exact sequence in the Leibniz 
theory maps to the ISB exact sequence in the cyclic theory.  
For an algebra $A$ over a characteristic zero field, this map 
between exact sequences acquires the form
$$ \CD
H^{\rm{rel}}_{n-1}(A) @>>>  HL_{n+1}(A) @>>> H^{\rm{Lie}}_{n+1}(A)
     @>>> H^{\rm{rel}}_{n-2}(A)  @>>> \\
@VVV \ @VV{\varphi_*}V \ @VV{\theta_*}V \ @VVV 	  \\
HC_{n-1}(A) @>>B> HH_n (A) @>>I> HC_n (A) @>>S> HC_{n-2}(A) @>>B>
\endCD  $$
where $H^{\rm{Lie}}_*(A)$ denotes Lie-algebra homology, and 
$H^{\rm{rel}}_*(A)$ is the homology of a mapping cone.  Every map
in the above diagram can be identified in the case of $gl(A)$.  

In the second half of the paper we study the above diagram in the
case of the universal enveloping algebra of formal vector fields.  
Let $W^P_1$ denote the Lie algebra of polynomial vector fields in 
one variable, and let $U(W^P_1)$ be its enveloping algebra.
The Godbillon-Vey invariant arises as a non-zero class in 
$HC_* (U(W^P_1))$.  We show that this class is detected by the map
$$ I \circ \varphi_* : HL_{*+1}(U(W^P_1)) \to HC_*(U(W^P_1)) .$$
The natural inclusion $W^P_1 \hookrightarrow U(W^P_1)$ induces
a map on Leibniz homology
$$ HL_*(W^P_1) \to HL_*(U(W^P_1))  , $$
and the composition
$$ HL_{*+1}(W^P_1) \to HC_{*}(U(W^P_1)) $$
is shown to detect the Godbillon-Vey invariant as well.

\section{A Map Between Exact Sequences}

To establish our notations and sign conventions, we begin with a
few definitions.  Let $\g$ be a Lie algebra over a commutative ring
$k$, or more generally consider $\g$ as a Leibniz algebra \cite{LodayP}.
The Leibniz homology of $\g$, written $HL_* (\g )$, is the homology
of the chain complex $CL_*(\g )$:
$$ \CD
k @<0<< \g @< [\ ,\ ]<< \g^{\ot 2} @<<< \ldots @<<< \g^{\ot (n-1)} @<d<< 
\g^{\ot n} @<<< \ldots \, ,
\endCD  $$
where, for $(g_1, \, g_2, \, \ldots , \, g_n ) \in \g^{\ot n}$, 
\begin{equation} \label{1}
\begin{split}
& d(g_1, \, g_2, \, \ldots , \, g_n) = \\
& \sum_{1 \leq i < j \leq n} (-1)^j \, (g_1, \, g_2, \, \ldots, \,
g_{i-1}, \, [g_i, \, g_j], \, g_{i+1},\, \ldots , \, \hat{g}_j, \, \ldots, \, g_n).
\end{split}
\end{equation}
Let $\Lambda^* (\g )$ be the chain complex for the Lie-algebra homology
of $\g$:
$$ \CD
k @<0<< \g @<[\ ,\ ]<< \g^{\wedge 2} @<<< \ldots @<<< \g^{\wedge (n-1)} @<d<<
\g^{\wedge n} @<<< \ldots \, , 
\endCD  $$
where $d$ is given essentially by equation \eqref{1}, with exterior powers 
in place of tensor powers.  The canonical projection
$$ \pi_1 : \g^{\ot n} \to \g^{\wedge n} $$
is a map of chain complexes and yields a short exact sequence
\begin{equation} \label{2}
\CD
0 @>>> ({\rm{ker}}\, \pi_1)_* [2] @>>> CL_*(\g ) @>{\pi_1}>>
\Lambda^*(\g ) @>>> 0 
\endCD
\end{equation}
where $({\rm{ker}}\, \pi_1)_n [2] = {\rm{ker}}[ \pi_1: \g^{\ot (n+2)} \to
\g^{\wedge (n+2)}]$.  Letting $H_*^{\rm{rel}}(\g )$ denote the homology
of the complex $({\rm{ker}} \, \pi_1 )_* [2]$, Pirashvili \cite{Pirashvili}
uses \eqref{2} to study the long exact sequence
$$ \CD
\cdots @>>> H^{\rm{rel}}_{n-2}(\g ) @>>> HL_n(\g ) @>>> H^{\rm{Lie}}_n(\g )
@>\partial>> H^{\rm{rel}}_{n-3}(\g) \\
\cdots @>>> H^{\rm{rel}}_0(\g ) @>>> HL_2(\g ) @>>> H^{\rm{Lie}}_2(\g ) 
@>>> 0 \, .
\endCD $$

In the sequel $k$ is a characteristic zero field, and $A$ a unital
algebra over $k$.  Consider $A$ as a Lie algebra via
$$ [a, \ b] = ab - ba, \ \ \ a, \ b \in A  . $$
Let $CHH_n(A) = (A^{\ot (n+1)}, \ b)$, $n \geq 0$, be the chain
complex for the Hochschild homology of $A$, written $HH_*(A)$, and
let
$$  C^{\lambda}_n(A) = (A^{\ot (n+1)}/{\scriptstyle{\bz/(n+1)}}, \ b), 
\ \ \ n \geq 0 , $$
be the chain complex for the cyclic homology of $A$ in the characteristic
zero setting.  See \cite[1.1.3, 2.1.4]{LodayC} for the details of these
constructions.  The projection 
$$ \pi_2 : A^{\ot (n+1)} \to A^{\ot (n+1)}/{\scriptstyle{\bz/(n+1)}} $$
is a map of chain complexes and yields a short exact sequence
\begin{equation} \label{3}
\CD
0 @>>> ({\rm{ker}}\, \pi_2)_* @>>> CHH_*(A) @>{\pi_2}>> C^{\lambda}_*(A)
@>>> 0 \, .
\endCD
\end{equation}
It is known that $H_n( ({\rm{ker}} \, \pi_2)_* ) \simeq
HC_{n-1}(A)$, and sequence \eqref{3} yields the long exact
sequence
$$ \CD
HC_{n-1}(A) @>B>> HH_n(A) @>I>> HC_n(A) @>S>> HC_{n-2}(A) @>B>>
\endCD $$
originally discovered by Connes \cite{Connes}.  Above, $I = (\pi_2)_*$. 

There is a chain map from the Lie-algebra homology complex to the 
cyclic complex \cite{LodayQ} \cite[10.2.3]{LodayC}  
\begin{equation} 
\begin{split}
  & \theta : \Lambda^{*+1} (A) \to C^{\lambda}_* (A)  \\
  & \theta(a_0 \wedge a_1 \wedge \, \ldots \, \wedge a_n) = \sum_{\sigma \in
  S_n} ({\rm{sgn}} \, \sigma) (a_0, \, a_{\sigma^{-1}(1)}, \,
  a_{\sg^{-1}(2)}\, \ldots , \, a_{\sigma^{-1}(n)}),
\end{split}
\end{equation}
which induces a homomorphism
$$  \theta_* : H^{\rm{Lie}}_{*+1}(A) \to HC_*(A)  . $$
There is a similar chain map from the Leibniz to the Hochschild chain complex
\cite{Lodder}
\begin{equation} 
\begin{split}
  & \varphi : CL_{*+1} (A) \to CHH_* (A)  \\
  & \varphi(a_0, \, a_1, \, \ldots \,  a_n) = \sum_{\sigma \in
  S_n} ({\rm{sgn}} \, \sigma) (a_0, \, a_{\sigma^{-1}(1)}, \,
  a_{\sg^{-1}(2)}\, \ldots , \, a_{\sigma^{-1}(n)}),
\end{split}
\end{equation}
which induces the map
$$ \varphi_* :HL_{*+1} (A) \to HH_*(A) . $$
\begin{proposition} \label{4}
For an algebra $A$ over a characteristic zero field, there
is a natural map of long exact sequences
$$ \CD
H^{\rm{rel}}_{n-1}(A) @>>>  HL_{n+1}(A) @>(\pi_1)_*>> H^{\rm{Lie}}_{n+1}(A)
     @>>> H^{\rm{rel}}_{n-2}(A)  @>>>  \\
@VV{\al_*}V \ @VV{\varphi_*}V \ @VV{\theta_*}V \ @VV{\al_*}V   \\
HC_{n-1}(A) @>>B> HH_n (A) @>>I> HC_n (A) @>>S> HC_{n-2}(A) @>>B>
\endCD  $$
\end{proposition}
\begin{proof}
The proof follows from the map of short exact sequences
$$ \CD
0 @>>> ({\rm{ker}}\, \pi_1)_{*+1}[2] @>>> CL_{*+1}(A) @>\pi_1>>
\Lambda^{*+1}(A) @>>> 0  \\
@. @VV{\al}V @VV{\varphi}V @VV{\theta}V  \\
0 @>>> ({\rm{ker}}\, \pi_2)_{*} @>>> CHH_*(A) @>\pi_2>>
C^{\lambda}_{*}(A) @>>> 0 \, . 
\endCD  $$
It is easily checked that $\theta \circ \pi_1 = \pi_2 \circ \varphi$.  
The chain map 
$$ \al : ({\rm{ker}}\, \pi_1)_{*+1}[2] \to
({\rm{ker}}\, \pi_2)_* $$ 
is given by the restriction of $\varphi$ to ker$\, \pi_1$.  
\end{proof}

Recall that $\bz /(n+1)$ acts on $A^{\ot (n+1)}$ via
$$ t(a_0, \, a_1, \, \ldots , \, a_n) =
(-1)^n \, (a_n, \, a_0, \, a_1, \, \ldots, \, a_{n-1}) $$
and $b' : A^{\ot (n+1)} \to A^{\ot n}$ is given by
$$ b' (a_0, \, a_1, \, \ldots , \, a_n) =
\sum_{i=0}^{n-1} (-1)^i \, (a_0, \, a_1, \, \ldots, \, 
a_i \cdot a_{i+1}, \, \ldots, \, a_n). $$
For completeness, we state the following result.
\begin{lemma}
The map $\al_* : H^{\rm{rel}}_{n-1}(A) \to HC_{n-1}(A)$ is induced on
the chain level by $\frac{1}{n}b' h \varphi$, where
$$ h = \frac{-1}{n+1} \big( t + 2t^2 + 3t^3 + \cdots + nt^n \big). $$
\end{lemma}
\begin{proof}
Note that
\begin{equation*}
\begin{split}
& {\rm{ker}}\, \pi_2 : A^{\ot (n+1)} \to A^{\ot (n+1)}/ 
\langle 1-t \rangle \\
& \simeq {\rm{Im}} \, (1-t): A^{\ot (n+1)} \to A^{\ot (n+1)}.
\end{split} 
\end{equation*}
Over a characteristic zero field, there is an isomorphism of complexes
$$ h: ({\rm{Im}}(1-t), \ b) \to (A^{\ot (*+1)}/N , \ b'), $$
where $N = \sum_{i=0}^n t^i$.   Furthermore, 
$H_* (A^{\ot (*+1)}/N , \ b')$ may be computed from a double complex
which begins with the acyclic column $(A^{\ot (*+1)}, \ b')$ and 
contains the resolution of $\bz /(n+1)$ acting on $A^{\ot (n+1)}$ along
the rows.  The lemma now follows from \cite[2.1.5]{LodayC}.  
\end{proof}

Let
$$ gl(A) = \underset{\overset{\longrightarrow}{n}}{\rm{lim}} \, gl_n(A) $$ 
be the Lie algebra of infinite matrices over $A$ with finitely many
nonzero entries.  Recall the definition of the trace map
$$ {\rm{tr}}: gl(A)^{\ot (n+1)} \to A^{\ot (n+1)} $$
given in \cite[1.2.1]{LodayC}.  By Morita invariance, both maps
\begin{equation*}
\begin{split}
& {\rm{tr}}_* : HH_*(gl(A)) \to HH_*(A)  \\
& {\rm{tr}}_* : HC_*(gl(A)) \to HC_*(A)
\end{split}
\end{equation*}
are isomorphisms.  From \cite{LodayQ} and \cite{Cuvier} we have respectively
\begin{equation*}
\begin{split}
& H^{\rm{Lie}}_* (gl(A)) \simeq \Lambda(HC_*(A)[1])  \\
& HL_*(gl(A)) \simeq T(HH_*(A))[1]) 
\end{split}
\end{equation*}
and from \cite[10.2.19]{LodayC} \cite{Lodder} the maps 
\begin{equation*}
\begin{split}
& {\rm{tr}}_* \circ \theta_* : H^{\rm{Lie}}_{*+1}(gl(A)) \to HC_*(A)  \\
& {\rm{tr}}_* \circ \varphi_* : HL_{*+1}(gl(A)) \to HH_*(A)
\end{split}
\end{equation*}
are epimorphisms.  
\begin{proposition}
For $A$ a unital algebra over a characteristic zero field, the map
$$ {\rm{tr}}_* \circ \al_* :H^{\rm{rel}}_{*-1}(gl(A)) \to HC_{*-1}(A) $$
is surjective.
\end{proposition}
\begin{proof}
In the calculation of $H^{\rm{Lie}}_*(gl(A))$, after taking $gl(k)$
invariants on the chain level for Lie-algebra homology, the subcomplex
of primitives has homology isomorphic to $HC_*(A)[1]$ with 
isomorphism given by ${\rm{tr}}_* \circ \theta_*$ \cite[10.2.19]{LodayC}.
Also, in the calculation of $HL_*(gl(A))$, after taking $gl(k)$ 
invariants on the chain level for Leibniz homology, a similar subcomplex
has homology isomorphic to $HH_*(A)[1]$ with isomorphism given by
${\rm{tr}}_* \circ \varphi_*$ \cite{Lodder}.  The result now follows
from Proposition \eqref{4} and the 5-lemma.
\end{proof}

\section{The Godbillon-Vey Invariant}

In this section we show that the Godbillon-Vey invariant, as
detected by the Leibniz homology of formal vector fields, maps
naturally to the Godbillon-Vey invariant as detected by the
cyclic homology of the universal enveloping algebra of these
vector fields.  Let $\langle \frac{d}{dx} \rangle$ be the real
one-dimensional vector space on the symbol $\frac{d}{dx}$, and
let 
$$ W^P_1 = \br[x] \underset{\br}{\ot} \langle
{\textstyle{\frac{d}{dx}}} \rangle $$
denote the Lie algebra of polynomial vector fields with bracket on
monomials given by 
$$ [x^n \, {\textstyle{\frac{d}{dx}}}, \ x^k \, {\textstyle{\frac{d}{dx}}}] 
= (k-n)\, x^{n+k-1} \, {\textstyle{\frac{d}{dx}}}. $$ 
The universal enveloping algebra $U(W^P_1)$ may also be considered as
a Lie algebra, and the inclusion
$$ W^P_1 \hookrightarrow U(W^P_1) $$ 
induces maps
\begin{equation} \label{5}
\begin{split}
& \mu_1 : HL_*(W^P_1) \to HL_*(U(W^P_1))  \\
& \mu_2 : H^{\rm{Lie}}_*(W^P_1) \to H^{\rm{Lie}}_*(U(W^P_1)) .
\end{split}
\end{equation}
We prove that the composition
$$ HL_*(W^P_1) \overset{\mu_1}{\longrightarrow} HL_*(U(W^P_1)) 
\overset{I \circ \varphi_*}{\longrightarrow} HC_{*-1}(U(W^P_1)) $$
detects the Godbillon-Vey invariant.

Recall that for a Lie algebra $\g$ over a characteristic zero field,
$HH_*(U(\g ))$ may be computed from the complex
$$ (S(\g ) \ot \Lambda^* (\g ), \ \delta) , $$
where $S( \g )$ is the symmetric algebra on $\g$ \cite[3.3.6]{LodayC}.
For $a \in S( \g )$, and $g_i \in \g$,
$$ \delta : S(\g)\ot \g^{\wedge n} \to S(\g )\ot \g^{\wedge (n-1)} $$
is given by
\begin{equation*}
\begin{split}
& \delta (a \ot g_1 \wedge g_2 \wedge \, \ldots \, \wedge g_n) =  \\
& \sum_{i=1}^n (-1)^{i+1} \{a, \, g_i \} \ot g_1 \wedge g_2 \wedge
\, \ldots \, \hat{g}_i \, \ldots \, \wedge g_n \ + \\
& \sum_{1 \leq i < j \leq n} (-1)^{j+1} a \ot g_1 \wedge \, \ldots
\, g_{i-1} \wedge [g_i, \, g_j] \wedge g_{i+1} \, \ldots \, \hat{g}_j
\, \ldots \, \wedge g_n \, ,
\end{split}
\end{equation*}
where $\{ a, \, g_i \}$ is the Poisson bracket.  The cyclic homology
groups $HC_* (U(\g ))$ may be computed from the double complex
$$ (S(\g ) \ot \Lambda^* (\g ), \ \delta, \ d) \simeq
( \Omega^*_{S(\g )}, \ \delta, \ d) , $$
which is often called a mixed complex \cite{Kassel} \cite[3.3.7]{LodayC}.
Above, $\Omega^*_{S(\g )}$ denotes the algebra of K\"ahler differentials with
differential operator $d$.  

We apply these results to $\g = W^P_1$ and work with a complex 
quasi-isomorphic to $(\Omega^*_{S(\g )}, \ \delta, \ d)$.  Let 
$$ e_i = x^{i+1}{\textstyle{\frac{d}{dx}}} \in W^P_1, \ \ \ 
i= -1,\ 0,\ 1,\ 2, \ \ldots \, .$$
Then $[e_i ,\, e_j] = (j-i)\, e_{i+j}$, and
$$  S(W^P_1) \simeq \br [e_{-1}, \, e_0, \, e_1, \, e_2, \, \ldots \, ] . $$
\begin{lemma}
The complex $(\Omega^*_{S(W^P_1)}, \ \delta, \ d)$ is quasi-isomorphic
to $(\Omega^*_{S(W^P_1)}/({\rm{Im}}\, d), \ \delta)$.
\end{lemma}
\begin{proof}
Arrange the double complex $(\Omega^*_{S(W^P_1)}, \ \delta, \ d)$ 
as
$$ \CD
@V{\delta}VV  @V{\delta}VV  \\
S(W^P_1) \ot (W^P_1)^{\wedge n} @<d<< 
S(W^P_1) \ot (W^P_1)^{\wedge (n-1)} @<d<<  \\
@V{\delta}VV  @V{\delta}VV  \\
S(W^P_1) \ot (W^P_1)^{\wedge (n-1)} @<d<< 
S(W^P_1) \ot (W^P_1)^{\wedge (n-2)} @<d<<  \\
@V{\delta}VV  @V{\delta}VV
\endCD $$
and filter this complex by rows.  The $E^0_{*, \, n}$ term of the
resulting spectral sequence becomes
$$ S(W^P_1) \ot (W^P_1)^{\wedge n} \overset{d}{\longleftarrow}
S(W^P_1) \ot (W^P_1)^{\wedge (n-1)} \overset{d}{\longleftarrow}
\ \cdots \, , $$
where $d$ is the exterior derivative on polynomials in the $e_i$'s.  
Let $V_k$ be the real vector space with basis
$$ e_{-1}, \ e_0, \ e_1, \ \ldots, \ e_k . $$
The homology of the complex
$$ S(V_k) \ot V^{\wedge n}_k \overset{d}{\longleftarrow}
S(V_k) \ot V^{\wedge (n-1)}_k \overset{d}{\longleftarrow}
\ \cdots  $$
becomes $S(V_k) \ot V^{\wedge n}_k /({\rm{Im}}\, d)$, since
for a polynomial algebra over a characteristic zero filed, 
every closed form is exact.  Note that
$$ S(W^P_1) = \underset{\overset{\longrightarrow}{k}}{\rm{lim}}\, S(V_k)$$  
and the $E^1$ term of the spectral sequence becomes
\begin{equation*}
\begin{split}
& E^1_{0, \, n} = S(W^P_1) \ot (W^P_1)^{\wedge n}/{\rm{Im}}\, d  \\
& E^1_{i, \, n} = 0, \ \ \ i \geq 1 .
\end{split}
\end{equation*}
The boundary map $E^1_{0,\, n} \to E^1_{0, \, n-1}$ is induced by
$\delta$.
\end{proof}

Define a weight, wt, on monomials of $S(W^P_1)$ by
$$ {\rm{wt}}(e^{k_1}_{i_1} \cdot e^{k_2}_{i_2} \cdot
\, \ldots \, \cdot e^{k_m}_{i_m}) = k_1 i_1 + k_2 i_2 + \, 
\cdots \, + k_m i_m . $$
Extend wt to monomials of $\Omega^*_{S(W^P_1)} \simeq
S(W^P_1) \ot \Lambda^*(W^P_1)$ by
$$ {\rm{wt}} (f \ot de_{j_1} \wedge de_{j_2} \wedge 
\, \ldots \, \wedge de_{j_n}) =
{\rm{wt}}(f) + j_1 + j_2 + \, \cdots \, + j_n , $$
where $f$ is a monomial in $S(W^P_1)$.  Then $\delta$ preserves weight,
since
$$ [e_i, \, e_j] = (j-i) \, e_{i+j} , $$
while $d$ preserves weight, since $d(e^k_i) = k \, e^{k-1}_i \ot de_i$.
Thus, the complex 
$$ (S(W^P_1) \ot (W^P_1)^{\wedge *}/{\rm{Im}} \, d, \ \delta ) $$  
splits as a direct sum of complexes graded by the weights
$q = -1$, 0, 1, 2, $\ldots$ \, .  We prove that the homology of
$$ (S(W^P_1) \ot (W^P_1)^{\wedge *}/{\rm{Im}} \, d, \ \delta ) $$  
is carried by the weight zero subcomplex.
\begin{lemma}
For $q \neq 0$, the weight $q$ subcomplex of 
$(S(W^P_1) \ot (W^P_1)^{\wedge *}/{\rm{Im}} \, d, \ \delta )$  
is contractible.
\end{lemma}
\begin{proof}
For $n \geq 0$ define a contracting chain homotopy
$$ \gamma_n : S(W^P_1) \ot \Lambda^n (W^P_1)/ {\rm{Im}}\, d \to
S(W^P_1) \ot \Lambda^{n+1}(W^P_1)/{\rm{Im}}\, d $$
by
\begin{equation*}
\begin{split}
& \gamma_n (f \ot de_{j_1}\wedge de_{j_2} \wedge \, \ldots \,
\wedge \, de_{j_n}) =  \\
& \Big( \frac{(-1)^{n+1}}{{\rm{wt}}(f) + j_1 + j_2 + \ldots + j_n}
\Big) \, (f \ot de_{j_1}\wedge de_{j_2} \wedge \, \ldots \,
\wedge \, de_{j_n} \wedge de_0) ,
\end{split}
\end{equation*}
where $f$ is a monomial in $S(W^P_1)$.  One verifies that
$$ \delta \gamma_n + \gamma_{n-1} \delta = {\bf{1}}, $$
using the property of the Poisson bracket
$$ \{ f, \, e_0 \} = {\rm{wt}}(f) \cdot f . $$
Also, $\gamma [{\rm{Im}} \, d] \subseteq {\rm{Im}} \, d$,
which shows that $\gamma$ is well-defined on the quotient
by Im$\, d$.  
\end{proof}

For 
$$ z = e^{k_1}_{i_1} \cdot e^{k_2}_{i_2} \cdot
\, \ldots \, \cdot e^{k_m}_{i_m} \ot  
de_{j_1}\wedge de_{j_2} \wedge \, \ldots \, \wedge \, de_{j_n} $$
define the length of $z$ by $L(z) = k_1 + k_2 + \cdots + k_m$.  
Then
$$ \delta : S(W^P_1) \ot \Lambda^n(W^P_1) \to
S(W^P_1) \ot \Lambda^{(n-1)}(W^P_1) $$
preserves length and induces a splitting of the complex
$$ (S(W^P_1) \ot (W^P_1)^{\wedge *}/{\rm{Im}} \, d, \ \delta ) $$  
according to length, which is independent of the weight 
splitting.  The subcomplex with $q = 0$ and $L = 1$ has vector
space generators
$$ e_0 , \ \ \ e_{-1} \ot e_1 , \ \ \ e_{-1} \ot e_0 \wedge e_1 $$
in dimensions $n = 0$, 1, 2 respectively, while this subcomplex
is zero for $n \geq 3$.  Note that
$$ \delta (e_{-1}\ot e_1) = 2e_0, \ \ \ 
\delta (e_{-1}\ot e_0 \wedge e_1) = 0 $$
in $S(W^P_1) \ot \Lambda^{*}(W^P_1)/{\rm{Im}} \, d$.  
Thus, $e_{-1}\ot e_0 \wedge e_1$ represents a non-zero element in 
$$ H_2 (S(W^P_1) \ot \Lambda^{*}(W^P_1)/{\rm{Im}} \, d,\ \delta), $$
while this homology class is represented by the cycle
$$ \al = (e_{-1} \ot e_0 \wedge e_1, \, e^2_0) $$
in the double complex $(S(W^P_1) \ot \Lambda^{*}(W^P_1),\ \delta ,\ d)$. 
To compare the class of of $\al$ in cyclic homology to the 
Godbillon-Vey invariant in Lie algebra or Leibniz homology, we 
must map $\al$ to its corresponding class in the complex
$C^{\lambda}(U(W^P_1))$.  

For an algebra $A$, let $(C_*(A), \, b, \, B)$ denote the $(b, \, B)$
double complex for cyclic homology, and let $(C_*(A), \, b, \, b')$
denote the $(b, \, b')$ double complex \cite[2.1]{LodayC}.  Recall that 
for a Lie algebra $\g$ over a characteristic zero field, there is a
chain map \cite[3.3.7]{LodayC}
$$ f: (S(\g ) \ot \Lambda^*(\g ), \ \delta, \ d) \to
(C_*(U(\g )), \ b, \ B) $$
which induces an isomorphism on homology.  Moreover, there is a 
quasi-isomorphism
$$ \iota (C_*(A), \ b, \ B) \hookrightarrow
(C_*(A), \ b, \ b') $$
sending $x$ to $x \oplus sNx$ \cite[2.1.7]{LodayC}, and a natural
projection
$$ p: (C_*(A), \ b, \ b') \to C^{\lambda}_*(A) $$
which is a quasi-isomorphism in characteristic zero.
\begin{theorem}
The composition
$$ (p \circ \iota \circ f) : (S(W^P_1)\ot \Lambda^*(W^P_1),
\, \delta, \, d) \to C^{\lambda}_*(U(W^P_1)) $$
sends $\al = (e_{-1}\ot e_0 \wedge e_1, \, e^2_0)$ to
$$ e_{-1}\ot e_0 \ot e_1 - e_{-1}\ot e_1 \ot e_0 + 1\ot e_0 \ot e_0
- 1 \ot 1 \ot e^2_0 \in U(W^P_1)^{\ot 3}/\langle 1-t \rangle . $$
Moreover, in $HC_2(U(W^P_1))$
$$ (p \circ \iota \circ f)_* (\al ) = {\textstyle{\frac{1}{3}}}\, 
\sum_{\sigma \in S_3} ({\rm{sgn}}\, \sigma) 
\, (e_{\sigma (-1)} \ot e_{\sigma (0)} \ot e_{\sigma (1)}), $$
where the latter element can be written as 
$2(e_{-1} \wedge e_0 \wedge e_1)$.
\end{theorem}
\begin{proof}
>From \cite[3.3.7]{LodayC} we have
$$ f(\al ) = 
(e_{-1}\ot e_0 \ot e_1 - e_{-1}\ot e_1 \ot e_0 + 1\ot e_0 \ot e_0
- 1 \ot 1 \ot e^2_0, \ e^2_0) . $$
Thus, the result for $(p \circ \iota \circ f)(\al )$ follows.  Note 
that in $C^{\lambda}_*(U(W^P_1))$,
$$ b( -(1 \ot 1 \ot \ot e_0 \ot e_0)) =
1 \ot e_0 \ot e_0 - 1 \ot 1 \ot e^2 , $$
using equivalences up to a cyclic shift.  Thus,
$$ (p \circ \iota \circ f)_* (\al ) = 2(e_{-1} \wedge e_0 \wedge e_1). $$
\end{proof}
\begin{corollary}
The composition
$$ H^{\rm{Lie}}_*(W^P_1) \overset{\mu_2}{\longrightarrow}
H^{\rm{Lie}}_*(U(W^P_1)) \overset{\theta_*}{\longrightarrow}
HC_{*-1}(U(W^P_1)) $$
detects the Godbillon-Vey invariant, i.e.,
$$ H^{\rm{Lie}}_3 (W^P_1) \to HC_2(U(W^P_1)) $$
is injective.
\end{corollary}
\begin{proof}
The Godbillon-Vey class is represented by the cycle
$e_{-1} \wedge e_0 \wedge e_1$ in $H^{\rm{Lie}}_*(W^P_1)$ 
\cite{Godbillon}.  From the definition of $\theta$, we have
$$ (\theta_* \circ \mu_2)\, (e_{-1} \wedge e_0 \wedge e_1) = 
   2(e_{-1} \wedge e_0 \wedge e_1) . $$
\end{proof}
\begin{corollary}
The composition
$$ HL_*(W^P_1) \overset{(\pi_1)_*}{\longrightarrow}
H^{\rm{Lie}}_*(W^P_1) \longrightarrow 
HC_{*-1}(U(W^P_1)) $$
detects the Godbillon-Vey invariant.
\end{corollary}
\begin{proof}
>From \cite{LodderA} it follows that 
$$ (\pi_1)_* : HL_3 (W^P_1) \to H^{\rm{Lie}}_3(W^P_1) $$
is an isomorphism.  A choice of a vector-space generator of 
$HL_3(W^P_1)$ is 
$$ \gamma = {\textstyle{\frac{1}{2}}} e_{-1} \ot e_0 \ot e_1
- {\textstyle{\frac{1}{2}}} e_{0} \ot e_{-1} \ot e_1
+ {\textstyle{\frac{1}{6}}} e_{-1} \ot e_{-1} \ot e_2 . $$
Then $(\pi_1)_* (\gamma ) = e_{-1} \wedge e_0 \wedge e_1$.
\end{proof}
\begin{corollary}
The Godbillon-Vey invariant is represented by a
non-zero class in each homology group of the commutative 
diagram
$$ \CD
HL_3(W^P_1) @>(\pi_1)_*>> H^{\rm{Lie}}_3 (W^P_1)  \\
@V{\mu_1}VV  @VV{\mu_2}V  \\
HL_3(U(W^P_1)) @>(\pi_1)_*>> H^{\rm{Lie}}_3(U(W^P_1)) \\
@V{\varphi_*}VV  @VV{\theta_*}V  \\
HH_2(U(W^P_1)) @>I>> HC_2(U(W^P_1)) .
\endCD $$
\end{corollary}

\end{document}